# Physics Informed Piecewise Linear Neural Networks for Process Optimization

Ece Serenat Koksal[a], Erdal Aydin[a,b,*]

[a]*Department of Chemical and Biological Engineering, Koç University, Istanbul 34450, Turkey*

[b]*Koç University TUPRAS Energy Center (KUTEM), Koç University, Istanbul, 34450, Turkey*

## ABSTRACT

Constructing first-principles models is usually a challenging and time-consuming task due to the complexity of the real-life processes. On the other hand, data-driven modeling, and in particular neural network models often suffer from issues such as overfitting and lack of useful and high-quality data. At the same time, embedding trained machine learning models directly into the optimization problems has become an effective and state-of-the-art approach for surrogate optimization, whose performance can be improved by physics-informed training. In this study, it is proposed to upgrade piece-wise linear neural network models with physics informed knowledge for optimization problems with neural network models embedded. In addition to using widely accepted and naturally piece-wise linear rectified linear unit (ReLU) activation functions, this study also suggests piece-wise linear approximations for the hyperbolic tangent activation function to widen the domain. Optimization of three case studies, a blending process, an industrial distillation column and a crude oil column are investigated. For all cases, physics-informed trained neural network based optimal results are closer to global optimality. Finally, associated CPU times for the optimization problems are much shorter than the standard optimization results.

---

* Corresponding author.
E-mail address: eaydin@ku.edu.tr (E. Aydin).

## 1. Introduction

Artificial neural networks are mathematical expressions which are able to efficiently generalize the input-output relations of a process, inspired by brain neuron connections. Typically, classical neural networks are fully data-driven black-box models in which the first-principles model of the process is absent. Despite the black-box nature, artificial neural networks have been widely used in literature in recent years since they are useful, especially for modeling of the complicated systems when data is present.

The complexity of chemical engineering systems gives birth to the need for data-driven models, especially artificial neural networks due to their capability to approximate the nonlinear nature. However, the black-box nature of machine learning applications, missing or low quality of data and problems such as overfitting creates the requirement for hybrid models. For instance, in chemical engineering applications such as batch and semi-batch production, data is limited and noisy due to the dynamic nature (Thebelt et al., 2022). Even in continuous processes, many quality related intermediate product specifications cannot be measured online, contributing to many issues related to online optimization and control.

In recent years, physics-informed hybrid models have become a solution to respond to the aforementioned issues. Stiasny et al. (2021) plugs implicit Runge-Kutta integration scheme directly into the physics-informed ANN training model of a dynamic power system. Nellikkath & Chatzivasileiadis (2021a) imports physics-informed constraints to the DC optimal flow power model for safety and minimizes worst-case violations. Another study by Nellikkath & Chatzivasileiadis (2021b) introduces AC power flow equations inside the neural network and reduces the worst-case violations. Masi et al. (2021) constructs a thermodynamics-based artificial neural network model which replaces constitutive laws by encoding two laws of thermodynamics in the network to predict material response. Pun et al. (2019) introduces interatomic bonding to machine learning potentials to improve transferability of these and achieve better atomistic simulations. Haghighat & Juanes (2021) developed a toolbox called SciANN in Python in which scientific calculations can be embedded into artificial neural networks. HybridML is another open-source toolbox where it is possible to construct hybrid models with differential equations. The toolbox reports a case study where the spread of COVID-19 in German federal states is modeled (Merkelbach et al., 2022).

For ANNs, hyperbolic tangent (tanh), rectified linear unit (ReLU) and sigmoid function are widely used as activation functions in the hidden layers. Rectified linear unit consists of two linear units with a discontinuity at the origin, resulting in a natural piece-wise linear property. The piecewise linear nature of ReLU makes the optimization suitable for mixed integer linear programming (MILP). Grimstad & Andersson (2019) shows MILP formulations of neural networks with ReLU activation functions. Tsay et al. (2021) proposes MILP formulations for ReLU based neural networks where model tightness and size are balanced. Katz et al. (2020a) combines deep learning models and multiparametric programming using ReLU as activation function. The formulation is solved as a MILP, resulting in a decrease in computational burden. Di Martino et al. (2022) reformulates ReLU trained neural networks as MILP for process optimization and the model is used for derivation of Pareto fronts between two different objectives. Katz et al. (2020b) integrates ReLU based neural networks with explicit MPC to control a solar field. Yang & Bequette (2021) uses ReLU based input convex neural networks to convert optimization-based control problems to convex ones. By convex formulation, they improve robustness and performance of the controller. Optimization and machine learning toolkit (OMLT) is developed by Ceccon et al. (2022) which uses pre-trained machine learning models to construct optimization formulations. reluMIP toolbox solves MILP optimization problem of ReLU trained ANN models using Gurobi (Lueg et al., 2021).

Hyperbolic tangent function (tanh) may provide better approximation capability for highly nonlinear processes, with a potentially higher possibility to lead to better training and test performances. For hyperbolic tangent activation functions, Schweidtmann & Mitsos (2019) proposes reduced space formulations and relaxed form of tanh, which reduces the dimensionality by McCormick based algorithms for mixed integer nonlinear global optimization and using their in-house solver MaiNGO. For piecewise linear approximations of nonlinear functions, D'Ambrosio et al. (2010) proposes three different formulations, one-dimensional, triangle and rectangle method, and discusses the trade-off between approximation quality and computational effort related to each method. Vielma (2015) suggests MILP formulations by introducing auxiliary variables and using a linear programming (LP) based technique. Sridhar et al. (2013) uses modified version of standard MIP formulations for piecewise linear functions with an indicator variable turned on. Sildir & Aydin (2022) suggest 3, 5 and 7- piece linear forms of tanh for convex training and they demonstrate that such approximations result in quite tight lower bounds and similar ANN performances both for training and test.

This study combines the physics-informed trained piecewise linear neural networks with process optimization, where piece-wise linear forms of the nonlinear activations are approximated before the optimization problem is formulated as an MILP. The main contributions of the proposed work are: (*i*) extensions for neural network based MILP optimization studies limited to ReLU activation function, (*ii*) performance comparison of standard and physics-informed neural networks for process optimization, (*iii*) application of physics informed trained models for optimization using piecewise linear approximations. In Section 2, the methodology and mathematics behind the training and optimization procedures are explained. Section 3 consist of three case scenarios where the advantage of physics-informed neural networks on optimization is investigated. In subsections, different training algorithms and optimization formulations are also investigated, including piecewise formulations for hyperbolic tangent. The concluding remarks are given in Section 4.

# 2. Methodology

## 2.1. Physics-informed ANN: Training

A feed forward fully connected artificial neural network can be expressed as:

$$\hat{y} = f_1(Af_2(Bu + C) + D) \quad (1)$$

where $f_1$ and $f_2$ are output and hidden layer activation functions, $A$ and $B$ are weight matrices, $C$ and $D$ are bias vectors, and $u$ and $\hat{y}$ are input and predicted output vectors, respectively.

The weights and biases are estimated through nonlinear optimization, called as training, where the objective function minimizes mean squared error (MSE) expressed below:

$$\begin{gathered} min_{A,B,C,D} \frac{1}{N}\sum_{i=1}^{N}(f_1(Af_2(Bu_i + C) + D) - y_i)^2 \\ s.t. \\ LW \leq A, B, C, D \leq UW \\ LH \leq f_2(Bu + C) \leq UH \end{gathered} \quad (2)$$

where $y_i$ is the i$^{th}$ data point, $N$ is the number of samples in the training, $LW$ and $UW$ are the lower and upper bound for the weight and bias terms, and $LH$ and $UH$ are lower and upper bounds for the outputs of the hidden layer, respectively.

To introduce a physics informed term in the training, loss function can be proposed as bi-objective, by minimizing error and physics term together, or additional constrains with proper bounds can be included. Bi-objective loss function of training can be written as:

$$min_{A,B,C,D}\frac{1}{N}\sum_{i=1}^{N}(f_1(Af_2(Bu_i+C)+D)-y_i)^2+(\alpha p)$$

$$s.t. \tag{3}$$

$$LW\leq A,B,C,D\leq UW$$

$$LH\leq f_2(Bu+C)\leq UH$$

where α and $p$ are scalar weight value and the physics function term, respectively. When $p$ is introduced as a constrain, training problem, sometimes addressed as physics-constrained training, is expressed as follows:

$$min_{A,B,C,D}\frac{1}{N}\sum_{i=1}^{N}(f_1(Af_2(Bu_i+C)+D)-y_i)^2$$

$$s.t. \tag{4}$$

$$\alpha p\leq U$$

$$LW\leq A,B,C,D\leq UW$$

$$LH\leq f_2(Bu+C)\leq UH$$

where $U$ represents some upper bound related to physics of the process.

Since additional constrains may increase the complexity of the nonlinear optimization, model can be trained in two parts, using warm-starting. Firstly, weight and bias parameters can be obtained without a constrain to find an initial guess, which is standard training. Then, in the second part, constrains are introduced and optimal weight and bias parameters can be found, as follows:

$$min_{A,B,C,D}z=\frac{1}{N}\sum_{i=1}^{N}(f_1(Af_2(Bu_i+C)+D)-y_i)^2 \tag{5a}$$

$$s.t.$$

$$LW \leq A, B, C, D \leq UW$$

$$LH \leq f_2(Bu + C) \leq UH$$

$$min_{A,B,C,D} z = \frac{1}{N}\sum_{i=1}^{N}(f_1(Af_2(Bu_i + C) + D) - y_i)^2$$

$$s.t. \tag{5b}$$

$$z \leq Z$$

$$\alpha p \leq U$$

where $Z$ is the upper bound for MSE, $z$.

## 2.2. Optimization with physics-informed trained neural networks

The objective of standard process optimization is to find the optimal inputs of the process that minimize a chosen output subject to the first-principles process model and physical constraints. Since the detailed and physics-based process models are expensive to build and requires detailed process knowledge, surrogate optimization is a viable alternative when process data is available. Accordingly, physics-informed neural network models can be embedded into the optimization problems as equality constraints, including the optimal values of the weight and bias parameters of the trained neural network model. This way, instead of standard machine learning models, physics-informed machine learning models, and in particular physics-informed neural network models, may bring about better performance for surrogate optimization.

Using the ReLU activation function for neural network models and embedding the trained neural network model into the optimization problem directly contributes to convex optimization problems after using Big-M formulations due to the nature of the ReLU. In this setting, the physics informed training and optimization steps can be easily performed using standard solvers and toolboxes (Abadi et al., 2016; Ceccon et al., 2022; Lueg et al., 2021). However, when the model is trained using the nonlinear hyperbolic tangent as the activation function, finding the optimal solution with the trained neural network model is challenging as the resulting formulation is a non-convex NLP problem. Karuppiah & Grossmann (2006) proposed a spatial branch and contract algorithm for solving non-convex NLP formulations, where piece-wise linear estimators are used to approximate the non-convex terms to obtain the lower bound of the global optimum. Schweidtmann & Mitsos (2019) proposed reduced space formulations with convex and concave envelopes for deterministic global optimization for dealing with this issue in a rigorous way.

In this study, to ensure convexity for the surrogate optimization problem, hyperbolic tangent function is approximated using piece-wise linear segments. Accordingly, the physics-informed neural network can be trained again using standard toolboxes and the original hyperbolic tangent activation function. After that, optimization is performed with the piece-wise approximation of the hyperbolic tangent. In fact, sigmoid and other types of activation functions can also be used in the same setting. Fortunately, the authors already showed that such approximations often result in tight relaxations and accurate representations (Sildir & Aydin, 2022). Given that -4, -1, +1 and +4 are the breakpoints for hyperbolic tangent function (tanh), 3 and 5-piece linear approximations of tanh are given as follows (D'Ambrosio et al., 2010):

$$\tanh(x)_{3\,pcs} \approx \begin{cases} 0.08x - 0.68 & if -4 \leq x \leq -1 \\ 0.76x & if -1 \leq x \leq +1 \\ 0.08x + 0.68 & if -1 \leq x \leq +4 \end{cases} \tag{6}$$

$$\tanh\,(x)_{5\,pcs} \approx \begin{cases} 0.018x - 0.93 & if -4 \leq x \leq -2 \\ 0.20x - 0.56 & if -2 \leq x \leq -1 \\ 0.76x & if -1 \leq x \leq +1 \\ 0.20x + 0.56 & if +1 \leq x \leq +2 \\ 0.018x + 0.93 & if +2 \leq x \leq +4 \end{cases} \tag{7}$$

where $\tanh(x)_{3\,pcs}$ and $\tanh(x)_{5\,pcs}$ indicate 3-piece and 5-piece approximation to tanh, respectively. Fig. 1. shows 3 and 5-piece approximations to tanh. 5-piece form sustains closer approximation since it has more breakpoints and linear segments. However, when the problem is formulated as an MILP, number of variables enhances with increasing number of linear segments.

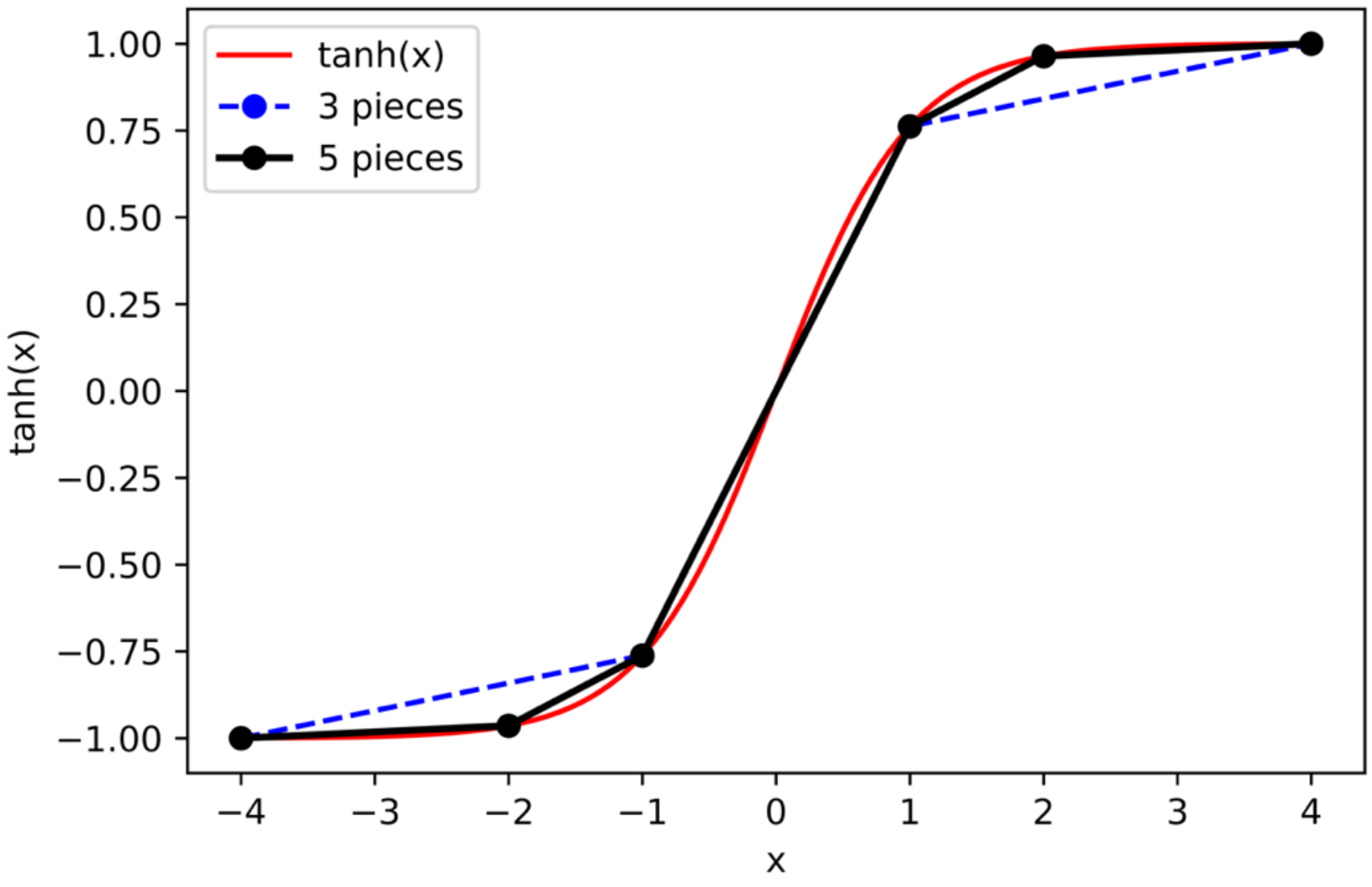

**Fig. 1.** Approximations to tanh.

Finally, the approximated tanh-based physics informed neural network model can be formulated into a mixed-integer linear optimization problem to find the optimal inputs that minimize the targeted output of the trained model using concepts in (Sildir & Aydin, 2022; Vielma, 2015) as follows:

$$\min_{u} y = f_1(A^* f_2(B^* u + C^*) + D^*)$$

$$s.t. \tag{8}$$

$$(-4\lambda_1) + (-\lambda_2) + (\lambda_3) + (4\lambda_4) = z_2$$

$$(-\lambda_1) + (-0.76\lambda_2) + (0.76\lambda_3) + (\lambda_4) = f_2'(z_2)$$

$$\lambda_1 + \lambda_2 + \lambda_3 + \lambda_4 = 1$$

$$\gamma_1 + \gamma_2 + \gamma_3 = 1$$

$$\delta_1 + \delta_2 + \delta_3 + \delta_4 = 1$$

$$(-4\delta_1) + (-\delta_2) + (\delta_3) + (4\delta_4) = z_1$$

$$(-\delta_1) + (-0.76\delta_2) + (0.76\delta_3) + (\delta_4) = f_1'(z_1)$$

$$-1 \leq f_1'(z_1) \leq 1$$

$$-1 \leq f_2'(z_2) \leq 1$$

$$-4 \leq z_1 \leq 4$$

$$-4 \leq z_2 \leq 4$$

$$\lambda_1 \leq \gamma_1$$

$$\lambda_2 \leq \gamma_1 + \gamma_2$$

$$\lambda_3 \leq \gamma_2 + \gamma_3$$

$$\lambda_4 \leq \gamma_3$$

$$\delta_1 \leq \beta_1$$

$$\delta_2 \leq \beta_1 + \beta_2$$

$$\delta_3 \leq \beta_2 + \beta_3$$

$$\delta_4 \leq \beta_3$$

$$\lambda_k \geq 0 \quad \forall k \in \{1, \ldots 4\}$$

$$\delta_k \geq 0 \quad \forall k \in \{1, \ldots 4\}$$

$$0 \leq \gamma_j \leq 1 \quad \forall j \in \{1,2,3\}$$

$$0 \leq \beta_j \leq 1 \quad \forall j \in \{1,2,3\}$$

$$\gamma_j \in \mathbb{Z} \quad \forall j \in \{1,2,3\}$$

$$\beta_j \in \mathbb{Z} \quad \forall j \in \{1,2,3\}$$

where $\lambda_j$ and $\delta_j$ are nonnegative auxiliary variables, $k$ is the number of nonnegative auxiliary variables for a single layer, $\gamma_j$ and $\beta_j$ are binary variables associated with breakpoints, $j$ is the number of linear segments of the approximation, $z_1$ and $z_2$ are the inputs for output and hidden layers which are defined between -4 and +4, $f_1'$ and $f_2'$ are tanh approximation function for output and hidden layers which are defined between -1 and +1, and $A^*$, $B^*$, $C^*$ and $D^*$ are the weight and bias parameters of the trained model, respectively.

Finally, the problem given by Eq. 8 can also be formulated by special ordered set Type-II (sos2) variables as follows (Sridhar et al., 2013):

$$
\begin{gathered}
\min_u y = f_1(A^* f_2(B^* u + C^*) + D^*) \\
s.t. \\
\lambda_1 + \lambda_2 + \lambda_3 + \lambda_4 = 1 \\
(-4\lambda_1) + (-\lambda_2) + (\lambda_3) + (4\lambda_4) = z_1 \\
(-\lambda_1) + (-0.76\lambda_2) + (0.76\lambda_3) + (\lambda_4) = f_1'(z_1) \\
\delta_1 + \delta_2 + \delta_3 + \delta_4 = 1 \\
(-4\delta_1) + (-1\delta_2) + (1\delta_3) + (4\delta_4) = z_2 \\
(-1\delta_1) + (-0.76\delta_2) + (0.76\delta_3) + (1\delta_4) = f_2'(z_2) \\
0 \le \lambda_k \le 1, k = 1,2,3,4 \\
0 \le \delta_k \le 1, k = 1,2,3,4 \\
-1 \le f_1'(z_1) \le 1 \\
-1 \le f_2'(z_2) \le 1 \\
-4 \le z_1 \le 4 \\
-4 \le z_2 \le 4 \\
\lambda \; and \; \delta \; are \; SOS2 \; variables.
\end{gathered}
\tag{9}
$$

The methodology is visualized in Fig. 2, where $u^*$ represent the optimal input set.

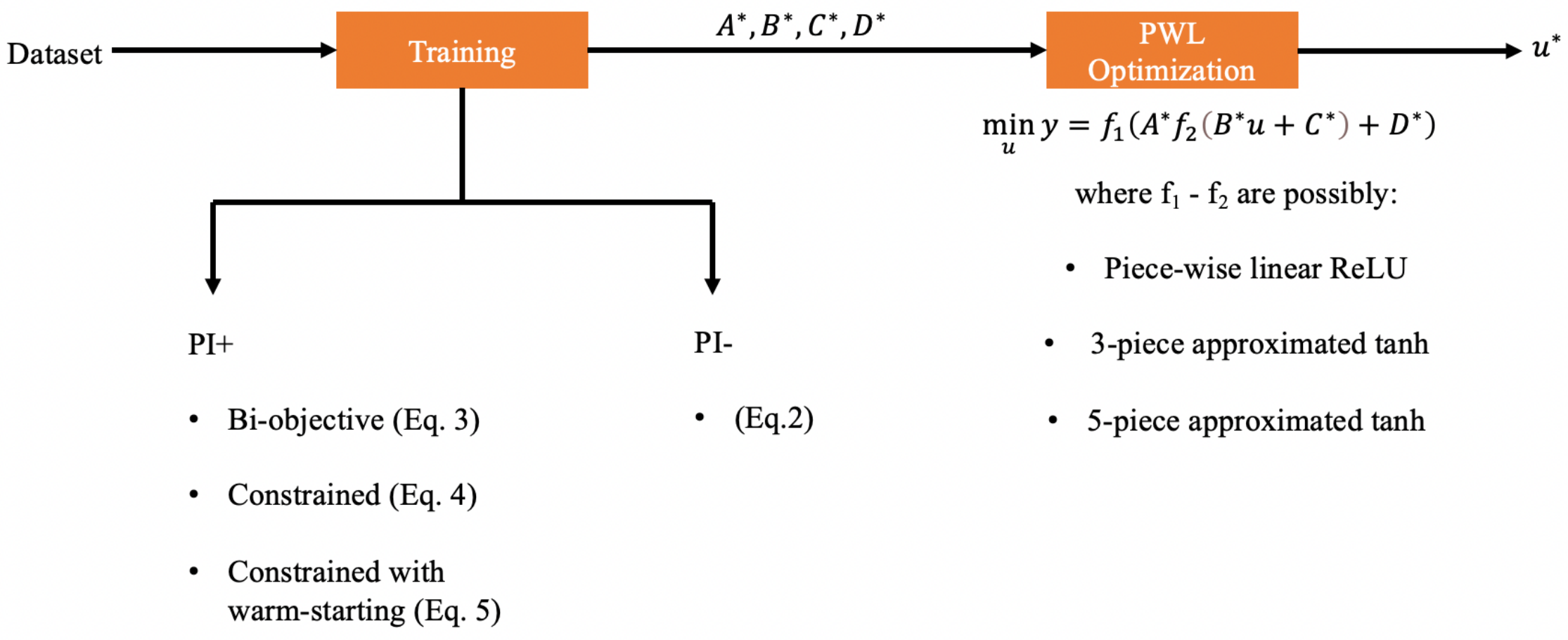

**Fig. 2.** Flowchart of the methodology.

## 3. Results and discussion

### 3.1. Case study 1: A simple blending process

A blending process of two streams in which two miscible liquids, are mixed, is shown in Fig. 3. Here, streams consist of a mixture of two components, A and B, and no reaction occurs. The inputs, $x_1$, $x_2$, $w_1$, and $w_2$ are the mole fraction of component A in stream 1, mol fraction of component A in stream 2, the molar flowrate of stream 1, and molar flowrate of stream 2, respectively. The outputs are represented by $x$ and $w$, mol fraction of component A and total flowrate in the exit stream, respectively. The input and output vectors for the ANN are given below:

$$u = \begin{bmatrix} x_1 \\ x_2 \\ w_1 \\ w_2 \end{bmatrix}$$

$$y = \begin{bmatrix} x \\ w \end{bmatrix}$$

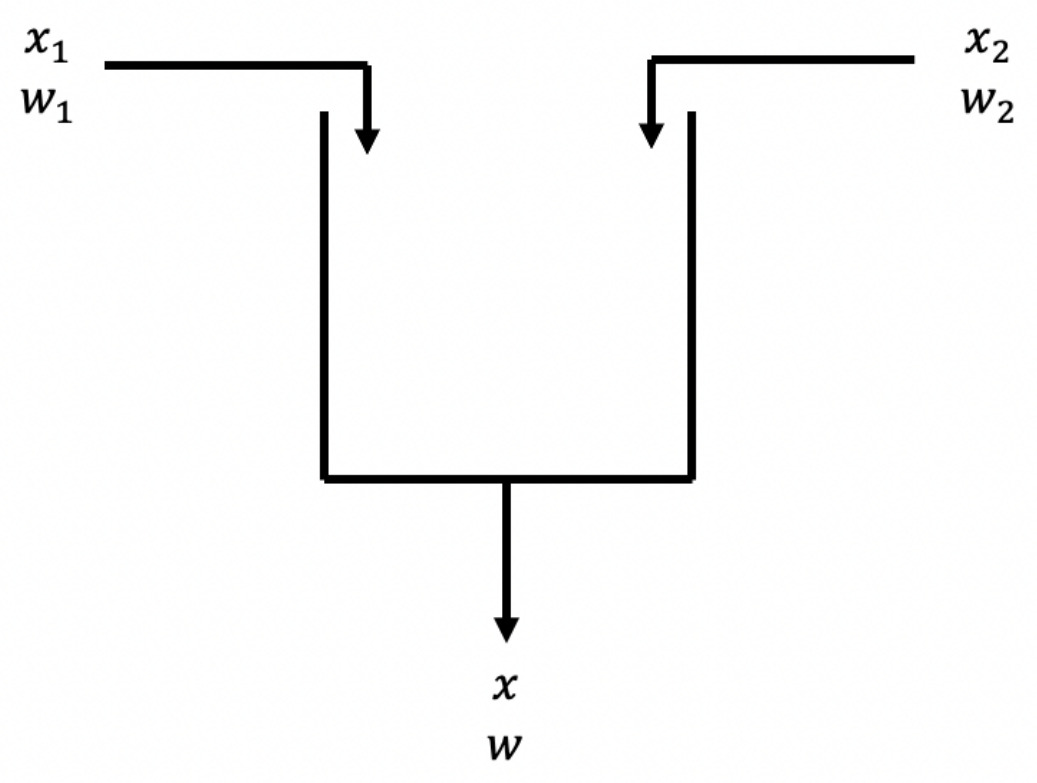

**Fig. 3.** Blending process.

The first-principles model equations of the process are given below:

$$\begin{gathered} x_1 w_1 + x_2 w_2 = xw \\ w_1 + w_2 = w \\ 0.1 \le x_1 \le 0.5 \\ 0.5 \le x_2 \le 1 \\ 0 \le w_1, w_2 \le 50 \end{gathered} \tag{10}$$

A synthetic dataset of size 100 is created by sampling the input data in Matlab using "rand" command based on the model in (10). The data is normalized between -1 and 1, white noise is added by "awgn" function and data is shuffled by "randperm". 80% of data is used for training. A single hidden layer with 5 neurons, in addition to the output layer, is used and tanh is chosen as the nonlinear activation function. The training procedures are as follows:

(i) Physics-uninformed (PI-) case, where the objective is to minimize MSE,

(ii) Physics-informed (PI+) case, where a physics term is embedded to the PI+ ANN by using the first term in Eq. 10, which is the component mass balance of A. The calculation of this term is given in Eq. 11:

$$p = \sum_{i=1}^{N} \left(\hat{y}_{1,i}\hat{y}_{2,i} - u_{i,1}u_{i,3} - u_{i,2}u_{i,4}\right)^2 \tag{11}$$

For physics uninformed (PI-) and physics-informed (PI+) cases, Eq. 2 and Eq.3 are chosen as objective functions, respectively. For the optimization part, the objective is minimizing the mol fraction of the outlet stream, given that the minimum flowrate of streams 1 and 2 must be 25 units. Considering tanh(x) is between -1 and +1 when x is defined between -4 and +4, the surrogate optimization problem with the exact form of tanh(x) can be formulated as given in Eq. 12:

$$\begin{gathered} \min_u y_1 = f_1(A^* f_2(B^* u + C^*) + D^*) \\ s.t. \\ -1 \leq f_1'(z_1) \leq 1 \\ -1 \leq f_2'(z_2) \leq 1 \\ 25 \leq u_3 \end{gathered} \tag{12}$$

$$25 \leq u_4$$

$$-4 \leq z_1 \leq 4$$

$$-4 \leq z_2 \leq 4$$

The optimization is performed for two different cases:

(i) The optimization is formulated as an NLP as shown in Eq. 12, which uses exact tanh for activation function, and BARON is used as the nonconvex, global solver. This is the surrogate optimization with the exact form of tanh and is expressed as $NLP_{tanh}$.

(ii) The problem in Eq. 12 is converted into an MILP, by using the 3 and 5 piece approximations to tanh as the activation function of the trained neural network, and solved by CPLEX. The formulations are expressed as $MILP_{3\text{-}pcs\text{-}tanh}$ and $MILP_{5\text{-}pcs\text{-}tanh}$, respectively. MILP formulation for 3-piece approximation to tanh for this case study is given in Eq. 13:

$$\begin{gathered}
\min_u y_1 = f_1(A^* f_2(B^* u + C^*) + D^*) \\
s.t. \\
\lambda_1 + \lambda_2 + \lambda_3 + \lambda_4 = 1 \\
(-4\lambda_1) + (-\lambda_2) + (\lambda_3) + (4\lambda_4) = z_1 \\
(-\lambda_1) + (-0.76\lambda_2) + (0.76\lambda_3) + (\lambda_4) = f_1'(z_1) \\
\delta_1 + \delta_2 + \delta_3 + \delta_4 = 1 \\
(-4\delta_1) + (-\delta_2) + (\delta_3) + (4\delta_4) = z_2 \\
(-\delta_1) + (-0.76\delta_2) + (0.76\delta_3) + (\delta_4) = f_2'(z_2) \\
0 \leq \lambda_k \leq 1, k = 1,2,3,4 \\
0 \leq \delta_k \leq 1, k = 1,2,3,4 \\
-1 \leq f_1'(z_1) \leq 1
\end{gathered} \tag{13}$$

$$-1 \leq f_2'(z_2) \leq 1$$

$$0 \leq u_3 \leq 1$$

$$0 \leq u_4 \leq 1$$

$$-4 \leq z_1 \leq 4$$

$$-4 \leq z_2 \leq 4$$

$$\lambda \ and \ \delta \ are \ SOS2 \ variables.$$

For PI- and PI+, 5 runs for each are reported with the training and test MSE, input set, objective and CPU time of the solver in Table 1. The reason for five trials is that ANN model is trained by local solvers. Due to nonconvexity of the training problem, training problem may converge to different local optima, so the trained networks may be different from each other. However, the differences in the optimal values, especially for the physics-informed cases should be small due to the increased reliability. Authors also suggested a convex-training method for training and feature selection for neural networks (Sildir & Aydin, 2022).

**Table 1.**

Training, test, and optimization performances for case study 1.

| *Training condition* | *Run #* | *Train MSE* | *Test MSE* | *Formulation* | $x_1$ | $x_2$ | $w_1$ | $w_2$ | $x$ | $w$ | *CPU time(s)* |
|---|---|---|---|---|---|---|---|---|---|---|---|
| PI- | 1 | 0.0062 | 0.0100 | $NLP_{tanh}$ | 0.12 | 0.5 | 35.9 | 25 | 0.28 | 70.0 | 10.2 |
| | | | | $MILP_{3\text{-}pcs\text{-}tanh}$ | 0.29 | 0.5 | 25 | 25 | 0.48 | 44.2 | 0.52 |
| | | | | $MILP_{5\text{-}pcs\text{-}tanh}$ | 0.26 | 0.5 | 47.9 | 25 | 0.30 | 62.9 | 0.88 |
| | 2 | 0.0022 | 0.0044 | $NLP_{tanh}$ | 0.1 | 0.5 | 49.6 | 25 | 0.24 | 78.9 | 1.50 |
| | | | | $MILP_{3\text{-}pcs\text{-}tanh}$ | 0.1 | 0.5 | 49.6 | 25 | 0.20 | 65.3 | 0.01 |
| | | | | $MILP_{5\text{-}pcs\text{-}tanh}$ | 0.1 | 0.5 | 49.6 | 25 | 0.18 | 65.3 | 0.03 |
| | 3 | 0.0032 | 0.0072 | $NLP_{tanh}$ | 0.1 | 0.5 | 38.5 | 25 | 0.29 | 68.5 | 7.15 |
| | | | | $MILP_{3\text{-}pcs\text{-}tanh}$ | 0.1 | 0.5 | 50.0 | 50 | 0.35 | 91.5 | 0.26 |
| | | | | $MILP_{5\text{-}pcs\text{-}tanh}$ | 0.1 | 0.5 | 42.4 | 25 | 0.27 | 79.3 | 0.53 |
| | 4 | 0.0033 | 0.0067 | $NLP_{tanh}$ | 0.1 | 0.5 | 50 | 25 | 0.22 | 82.7 | 2.53 |
| | | | | $MILP_{3\text{-}pcs\text{-}tanh}$ | 0.1 | 0.5 | 50 | 25 | 0.19 | 82.7 | 0.56 |
| | | | | $MILP_{5\text{-}pcs\text{-}tanh}$ | 0.1 | 0.5 | 50 | 25 | 0.16 | 82.9 | 0.36 |
| | 5 | 0.0036 | 0.0068 | $NLP_{tanh}$ | 0.1 | 0.5 | 50 | 25 | 0.29 | 81.0 | 6.71 |
| | | | | $MILP_{3\text{-}pcs\text{-}tanh}$ | 0.1 | 1 | 41.9 | 25 | 0.30 | 61.7 | 0.69 |
| | | | | $MILP_{5\text{-}pcs\text{-}tanh}$ | 0.1 | 0.57 | 25 | 25 | 0.28 | 52.2 | 0.49 |
| PI+ | 1 | 0.0034 | 0.0048 | $NLP_{tanh}$ | 0.1 | 0.5 | 50 | 25 | 0.24 | 79.9 | 1.07 |
| | | | | $MILP_{3\text{-}pcs\text{-}tanh}$ | 0.1 | 0.5 | 50 | 25 | 0.21 | 64.9 | 0.34 |
| | | | | $MILP_{5\text{-}pcs\text{-}tanh}$ | 0.1 | 0.5 | 50 | 25 | 0.21 | 65.4 | 0.54 |
| | 2 | 0.0060 | 0.0098 | $NLP_{tanh}$ | 0.1 | 0.5 | 50 | 25 | 0.21 | 78.9 | 1.61 |
| | | | | $MILP_{3\text{-}pcs\text{-}tanh}$ | 0.1 | 0.5 | 49.7 | 25 | 0.21 | 71.3 | 0.42 |
| | | | | $MILP_{5\text{-}pcs\text{-}tanh}$ | 0.1 | 0.5 | 49.7 | 25 | 0.21 | 72.8 | 0.35 |
| | 3 | 0.0030 | 0.0043 | $NLP_{tanh}$ | 0.1 | 0.5 | 50 | 25 | 0.23 | 79.9 | 1.72 |
| | | | | $MILP_{3\text{-}pcs\text{-}tanh}$ | 0.1 | 0.5 | 50 | 25 | 0.20 | 63.7 | 0.01 |
| | | | | $MILP_{5\text{-}pcs\text{-}tanh}$ | 0.1 | 0.5 | 50 | 25 | 0.21 | 78.0 | 0.50 |
| | 4 | 0.0064 | 0.0093 | $NLP_{tanh}$ | 0.1 | 0.5 | 50 | 25 | 0.26 | 77.8 | 0.75 |
| | | | | $MILP_{3\text{-}pcs\text{-}tanh}$ | 0.1 | 0.5 | 50 | 38.5 | 0.34 | 82.6 | 0.63 |
| | | | | $MILP_{5\text{-}pcs\text{-}tanh}$ | 0.1 | 0.5 | 50 | 38.5 | 0.34 | 82.6 | 0.26 |
| | 5 | 0.0030 | 0.0043 | $NLP_{tanh}$ | 0.1 | 0.5 | 50 | 25 | 0.23 | 79.9 | 1.61 |
| | | | | $MILP_{3\text{-}pcs\text{-}tanh}$ | 0.1 | 0.5 | 50 | 25 | 0.20 | 63.7 | 0.03 |
| | | | | $MILP_{5\text{-}pcs\text{-}tanh}$ | 0.1 | 0.5 | 50 | 25 | 0.21 | 78.0 | 0.89 |
| Global | | | | *NLP* | 0.1 | 0.5 | 50 | 25 | 0.23 | 75 | 0.01 |

Average training and test MSE values are 0.0037 and 0.0070 for PI- and 0.0043 and 0.0065 for PI+, respectively. In PI+ case, training performance decreases due to the noise in data, which causes the deviation from the mass conservation. As a result, introducing the physics term to the objective function decreases the training performance slightly. On the other hand, test performance increases since the physics term in training provides better prediction.

For the optimization part, in addition to surrogate optimization runs, the model is formulated as an NLP using the first-principles model and solved using the nonconvex global solver BARON. Globally optimal inputs and outputs are given at the last line of Table 1. The global solution implies that to achieve minimum mol fraction at the outlet stream, both the mol fractions and the flowrate of the concentrated stream, $x_1$, $x_2$, and $w_2$, should be minimum and the flowrate of the dilute stream, $w_1$, should be maximum, which can also be interpreted by engineering intuition.

When PI+ case is formulated as $NLP_{tanh}$, using the exact form of tanh for the surrogate optimization, and solved by the nonconvex global solver BARON, the values of the decision variables (inputs) always converge to the global optimum. On the other hand, optimal inputs of PI- trained neural networks with the exact tanh activation functions, $NLP_{tanh}$, may deviate from the global optimum, as in Run 1 and 3 of PI- case in Table 1. This shows that, physics-informed training itself is beneficial for surrogate optimization problems where neural network models are embedded.

Since physics informed term is based on component mass balance rather than overall mass conservation and ANN is trained with noisy data, the deviation is much more significant for the flowrate of the output stream, $w$, especially for the last 3 runs of PI- case. This is a promising result considering the quality and low volume of the data. In fact, when problem size is larger or when data quality is low, instead of increasing the number of neurons which leads to a more complex problem and possibly overfitting, physics-informed neural networks can be used to achieve more reliable surrogate optimization.

Another observation is that CPU times for $NLP_{tanh}$ formulated optimization of PI+ trained neural networks are shorter. Yet, please note that using exact forms of such highly nonlinear activation functions may not be practical for larger neural networks and deep learning (Schweidtmann & Mitsos, 2019). Accordingly, results of surrogate optimizations with PI+ piece-wise linear neural networks, contributing to MILP problems and being another contribution of the proposed work, are examined.

Results of PI- MILP formulated optimization, proposing approximations to tanh, deviates from the global optimum compared to the PI- formulated optimization where exact tanh is used ($NLP_{tanh}$). On the other hand, CPU times sharply decrease when optimization is formulated as an MILP, as expected. These differences can be significant even for small number of neurons. The formulation with exact tanh is nonconvex, which increases the CPU time especially when a global nonconvex solver is used. On the other hand, the piece-wise approximations to tanh converts the problem to an MILP. Due to the convexity, CPU time is much shorter.

This case study has small problem size and number of neurons, so the comparison of the global and MILP solutions is straightforward. For higher number of hidden layers (deep learning) and neurons, solution of the NLP problem is expected to become even more challenging and maybe impossible, whereas the same does not apply to the approximated MILP optimization problems with physics-informed trained neural networks (Schweidtmann & Mitsos, 2019). So, at the cost of small deviations from the global solution, MILP formulation with the approximated hyperbolic tangent activation functions gives much faster results for the investigated case study.

In addition to the computational advantages of MILP formulations, physics-informed training for approximated tanh activation function contributes to reduced gap between the surrogate optimization and exact global results. As given in Table 1, the deviations of MILP formulations from $NLP_{tanh}$ and global NLP are much more significant for physics uninformed, PI- cases. At the same time, for PI+, $MILP_{5\text{-pcs-tanh}}$ solutions are closer to $NLP_{tanh}$ than $MILP_{3\text{-pcs-tanh}}$, since 5-piece approximation has more linear segments and breakpoints, resulting in better approximation capacity. As a result, approximating the exact tanh with more pieces may sustain even negligible deviations from exact global optimum, at the cost of small increase in CPU time.

### 3.2. Case study 2: An industrial distillation column

For this case study, normalized data of the crude oil separation column in Tüpraş Kırıkkale Refinery is used for proprietary reasons. The crude oil (CR), liquified petroleum gas (LPG), light straight run naphtha (LSRN), heavy straight run naphtha (HSRN), light diesel (LD), heavy diesel (HD), kerosene (kero) and atmospheric residuum (RSD) flowrate data are present. Fig. 4 is a representation of the column.

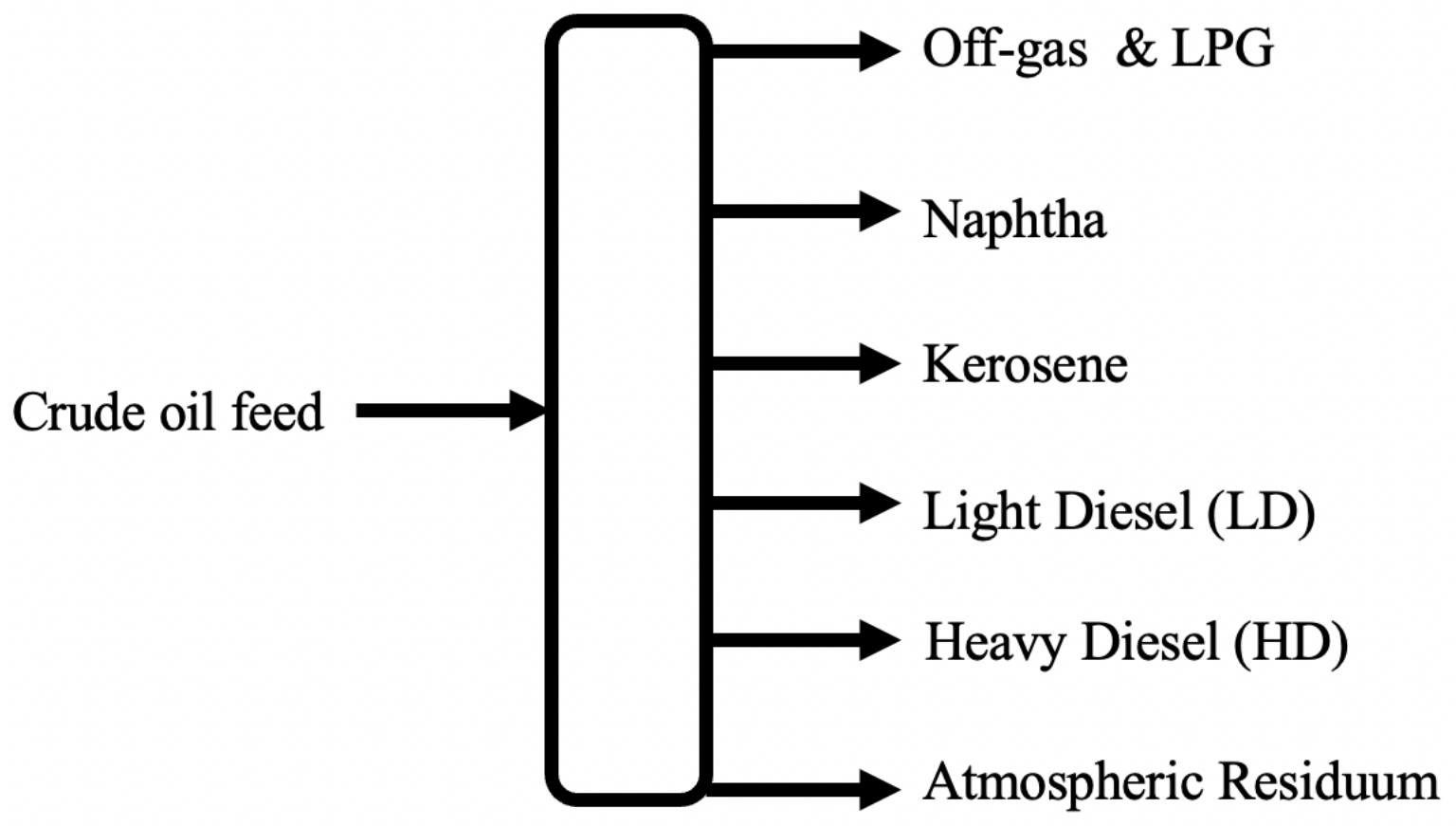

**Fig. 4.** Representation of the crude oil column.

ANN inputs are chosen as crude oil flowrate and other product flowrates, except residuum, which will be the predicted output of the ANN. The input and output vectors for the ANN become:

$$u = \begin{bmatrix} F_{CR} \\ F_{LPG} \\ F_{LSRN} \\ F_{HSRN} \\ F_{kero} \\ F_{LD} \\ F_{HD} \end{bmatrix}$$

$$y = [F_{RSD}]$$

where *F* represents the flowrate.

Since flowrate ranges are different, mass balance equation is imported by back normalization. Physics term, *p*, is defined as:

$$p = \frac{\sum_{i=1}^{N}\left(u_{1,i} - \sum_{m=2}^{7} u_{m,i} - \hat{y}_i\right)}{N} \tag{14}$$

where $\hat{y}$ is predicted value of flowrate of residuum. The data is normalized between -1 and 1 and the neural network is trained in two different ways:

(i) by minimizing the MSE as objective, as in Eq. 2, which is called the physics-uninformed (PI-) case,

(ii) by weighting MSE and $p$ terms and formulating the loss function bi-objective, as in Eq. 3, which is called as physics-informed (PI+) case.

The neural network has a single hidden layer with five neurons, the hyperbolic tangent is used as nonlinear activation function at the hidden and output layers. From 270 data points, 80% is used for training, and the rest is used for testing. Given that the crude oil flowrate at the inlet of the column is maximum, and considering tanh(x) values are between -1 and +1 when the decision variable interval is [-4, +4], the optimal inputs that minimizes the flowrate of the atmospheric residuum are found by the surrogate optimization problem given below:

$$\min_{u} y = f_1(A^* f_2(B^* u + C^*) + D^*)$$
$$s.t.$$
$$-1 \leq f_1'(z_1) \leq 1$$
$$-1 \leq f_2'(z_2) \leq 1 \tag{15}$$
$$u_1 = \eta_{CR,max}$$
$$-4 \leq z_1 \leq 4$$
$$-4 \leq z_2 \leq 4$$

where $\eta_{CR,max}$ is the maximum normalized crude oil value and is equal to 1 in this case study. Please note that detailed first-principles model of the investigated distillation column is unfortunately not present. However, since the column must be subject to mass conservation, assuming there is no leakage, a model is formulated based on a mass balance to obtain a prediction for the global solution and expressed as estimated global solution (*EGS*). The global solution based on the mass balance implies that, given that the crude oil is maximum, other products must also be maximum to obtain minimum atmospheric residuum.

The optimization problem in Eq. 15 is also formulated as an MILP by using 3 and 5 piece approximations for tanh. MILP formulations ($MILP_{3\text{-}pcs\text{-}tanh}$, $MILP_{5\text{-}pcs\text{-}tanh}$) and the NLP formulation with exact tanh in Eq. 15 ($NLP_{tanh}$) are solved by CPLEX and the nonconvex global solver BARON, respectively.

Due to the generic issue of nonconvexity in training, the model is trained five times and optimization is performed for these training instances separately to avoid the impacts of suboptimal training. Training and test performances are given in Table 2.

**Table 2.**

Training and test performance for case study 2.

| *Training condition* | *PI-* | | *PI+* | |
|---|---|---|---|---|
| Run # | Training MSE | Test MSE | Training MSE | Test MSE |
| 1 | $2.12\times10^{-5}$ | $5.38\times10^{-5}$ | $7.56\times10^{-6}$ | $2.44\times10^{-5}$ |
| 2 | $3.83\times10^{-5}$ | $1.23\times10^{-4}$ | $1.50\times10^{-5}$ | $4.58\times10^{-5}$ |
| 3 | $2.04\times10^{-5}$ | $4.90\times10^{-5}$ | $1.13\times10^{-5}$ | $3.47\times10^{-5}$ |
| 4 | $3.54\times10^{-4}$ | $6.28\times10^{-4}$ | $1.48\times10^{-5}$ | $4.94\times10^{-5}$ |
| 5 | $2.16\times10^{-5}$ | $7.12\times10^{-5}$ | $1.31\times10^{-5}$ | $5.46\times10^{-5}$ |
| Average | $9.11\times10^{-5}$ | $1.85\times10^{-4}$ | $1.23\times10^{-5}$ | $4.18\times10^{-5}$ |

For all runs in this case study, training and test performances improve when weighted physics term is added to the training objective function. For each trained model, results of three different optimization formulation are given in Table 3, where $\eta$ represents the normalized value of the variables.

**Table 3.**

Optimal inputs for case study 2.

| *Physics condition* | *Run* | *Formulation* | $\eta_{LPG}$ | $\eta_{LSRN}$ | $\eta_{HSRN}$ | $\eta_{Kero}$ | $\eta_{LD}$ | $\eta_{HD}$ | $\eta_{RSD}$ | *CPU time (s)* |
|---|---|---|---|---|---|---|---|---|---|---|
| PI- | 1 | $NLP_{tanh}$ | 1 | 1 | 1 | 1 | 1 | 1 | 0.21 | 0.76 |
| | | $MILP_{3\text{-}pcs\text{-}tanh}$ | 1 | 1 | 1 | 1 | 1 | 1 | 0.23 | 0.51 |
| | | $MILP_{5\text{-}pcs\text{-}tanh}$ | 1 | 1 | 1 | 1 | 1 | 1 | 0.12 | 0.53 |
| | 2 | $NLP_{tanh}$ | 0.42 | 1 | 1 | 1 | 1 | 1 | 0.17 | 22.0 |
| | | $MILP_{3\text{-}pcs\text{-}tanh}$ | 1 | 1 | 1 | 1 | 1 | 1 | 0.15 | 0.75 |
| | | $MILP_{5\text{-}pcs\text{-}tanh}$ | 1 | 1 | 1 | 1 | 1 | 1 | 0.03 | 0.61 |
| | 3 | $NLP_{tanh}$ | -0.20 | 1 | 1 | 1 | 1 | 1 | 0.25 | 9.58 |
| | | $MILP_{3\text{-}pcs\text{-}tanh}$ | -0.20 | 1 | 1 | 1 | 1 | 1 | -0.12 | 0.74 |
| | | $MILP_{5\text{-}pcs\text{-}tanh}$ | -0.20 | 1 | 1 | 1 | 1 | 1 | 0.10 | 0.56 |
| | 4 | $NLP_{tanh}$ | 1 | -0.10 | 1 | 1 | -0.04 | 1 | 0.53 | 69.5 |
| | | $MILP_{3\text{-}pcs\text{-}tanh}$ | 1 | -0.41 | 1 | 1 | 0.16 | 1 | 0.67 | 0.68 |
| | | $MILP_{5\text{-}pcs\text{-}tanh}$ | 1 | -0.29 | 1 | 1 | 0.08 | 1 | 0.66 | 0.571 |
| | 5 | $NLP_{tanh}$ | 1 | 0.65 | 1 | 1 | 1 | 1 | 0.25 | 1.71 |
| | | $MILP_{3\text{-}pcs\text{-}tanh}$ | 1 | 0.65 | 1 | 1 | 1 | 1 | 0.56 | 0.42 |
| | | $MILP_{5\text{-}pcs\text{-}tanh}$ | 1 | 0.65 | 1 | 1 | 1 | 1 | 0.22 | 0.56 |
| PI+ | 1 | $NLP_{tanh}$ | 1 | 1 | 1 | 1 | 1 | 1 | 0.19 | 1.53 |
| | | $MILP_{3\text{-}pcs\text{-}tanh}$ | 1 | -0.58 | 1 | 1 | 1 | 1 | -0.17 | 0.69 |
| | | $MILP_{5\text{-}pcs\text{-}tanh}$ | 1 | 1 | 1 | 1 | 1 | 1 | 0.20 | 0.52 |
| | 2 | $NLP_{tanh}$ | 1 | 1 | 1 | 1 | 1 | 1 | 0.20 | 4.91 |
| | | $MILP_{3\text{-}pcs\text{-}tanh}$ | -1 | 1 | 1 | 0.57 | 1 | 1 | 0.17 | 0.51 |
| | | $MILP_{5\text{-}pcs\text{-}tanh}$ | 1 | 1 | 1 | 1 | 1 | 1 | 0.17 | 0.44 |
| | 3 | $NLP_{tanh}$ | 1 | 1 | 1 | 1 | 1 | 1 | 0.20 | 3.25 |
| | | $MILP_{3\text{-}pcs\text{-}tanh}$ | -0.34 | 1 | 1 | 1 | 1 | 1 | 0.26 | 0.40 |
| | | $MILP_{5\text{-}pcs\text{-}tanh}$ | 1 | 1 | 1 | 1 | 1 | 1 | 0.18 | 0.51 |
| | 4 | $NLP_{tanh}$ | 1 | 1 | 1 | 1 | 1 | 1 | 0.20 | 1.26 |
| | | $MILP_{3\text{-}pcs\text{-}tanh}$ | -1 | 1 | 1 | 0.94 | 1 | 1 | 0.15 | 0.58 |
| | | $MILP_{5\text{-}pcs\text{-}tanh}$ | 1 | 1 | 1 | 1 | 1 | 1 | 0.17 | 0.61 |
| | 5 | $NLP_{tanh}$ | 1 | 1 | 1 | 1 | 1 | 1 | 0.20 | 1.37 |
| | | $MILP_{3\text{-}pcs\text{-}tanh}$ | 1 | 1 | 1 | 1 | 1 | 1 | 0.17 | 0.44 |
| | | $MILP_{5\text{-}pcs\text{-}tanh}$ | 1 | 1 | 1 | 1 | 1 | 1 | 0.22 | 0.48 |
| *EGS* | | | 1 | 1 | 1 | 1 | 1 | 1 | 0.20 | 0.00 |

When optimization is formulated with the exact tanh activation function in the hidden layer of the neural network, as $NLP_{tanh}$, PI+ trained model always converges to the estimated global optimum except the small deviation in $\eta_{RSD}$ at Run 1 of PI+, similar to the first case study. Results also show that PI- case can also possibly converge to the estimated global solution, as in Run 1. On the other hand, other runs for PI- show that there is no consistency, and it is more likely that PI- trained neural network models do not usually converge to the estimated global solution. Hence, it can be stated that the impact of training on the optimization is much higher in the physics-uninformed surrogate optimization case as PI+ trained neural networks sustains both close-to-optimum results and consistency.

In addition to better accuracy with PI+, although there are deviations due to approximation for the piece-wise linear networks, results of MILP formulations are still close to the $NLP_{tanh}$ solution with the physics-informed training. In fact, deviations of PI+ piece-wise linear networks are even negligible compared to the PI- piece-wise linear neural networks for surrogate optimization. Moreover, CPU times show that MILP formulation requires much shorter computational time, due to convexity, compared to $NLP_{tanh}$. In other words, approximations yield tight relaxations in the presence of PI+ training. Furthermore, $MILP_{5\text{-}pcs\text{-}tanh}$ results are closer to global optimum than $MILP_{3\text{-}pcs\text{-}tanh}$ with more linear segments and breakpoints. As a result, physics-informed training not only brings about closer global results for optimization, but also computational advantages for larger neural network-based MILP optimization, as also observed in this case study.

### 3.3. Case study 3: A crude oil unit

Crude oil unit (CDU) is the initial separation plant in a refinery where crude oil with long carbon chains is converted to white products, mainly LPG, naphtha, kerosene, diesel, and VGO. From the bottom, atmospheric residuum is obtained. This product can be used as asphalt or sent to other units for further cracking, such as vacuum distillation. For this case study, the model in Alhajri et al. (2008) is used with some modifications. The cut point temperature ranges used in the model are tabulated in Table 4.

**Table 4.**

Cut point temperature ranges.

| Product | Boiling Range Lower Bound (F) | Boiling Range Upper Bound (F) |
|---|---|---|
| **LPG** | -48 | -40 |
| **Naphtha** | 230 | 380 |
| **Kerosene** | 330 | 520 |
| **Diesel** | 420 | 630 |
| **VGO** | 620 | 1050 |

The cut points of the products are found using Eq. 16:

$$Cut_s = \sum_{k=0}^{4} a_k (TE_{CDU,s})^k \qquad (16)$$

where $s$ is the products of the CDU unit, except *RSD* (Residuum), $TE_{CDU,s}$ is the cut point temperature of product $s$ and $a_k$ are constants given in Table 5 (Alhajri et al., 2008). Since cut is cumulative, cut of *RSD* is equal to 100.

**Table 5.**

Constants for cut point equation.

| | |
|---|---|
| $a_0$ | 4.04 |
| $a_1$ | -0.047 |
| $a_2$ | $3.25\times10^{-4}$ |
| $a_3$ | $-2.84\times10^{-7}$ |
| $a_4$ | $8.15\times10^{-11}$ |

The flowrate of each product is found by:

$$F_{CDU,s} = F_{CDU} \times \left(\frac{Cut_s - Cut_{s-1}}{100}\right) \tag{17}$$

At different crude oil flowrate and cut temperatures, a synthetic dataset size of 1000 is created by Latin-hypercube sampling ("lhsdesign" function in Matlab), and product flowrates at these inputs are obtained using the model equations. The ANN and output variables are listed below:

$$u = \begin{bmatrix} F_{CDU} \\ TE_{LPG} \\ TE_{Naphtha} \\ TE_{Kerosene} \\ TE_{Diesel} \\ TE_{VGO} \end{bmatrix} \qquad y = \begin{bmatrix} F_{LPG} \\ F_{Naphtha} \\ F_{Kerosene} \\ F_{Diesel} \\ F_{VGO} \\ F_{RSD} \end{bmatrix}$$

Different training algorithms are performed and their MSE values are given in Table 6.

**Table 6.**

Training and test performance for case study 3.

| *# of neurons in the hidden layer* | *Training algorithm* | *Hidden layer activation function* | *Training condition* | *Training MSE* | *Test MSE* |
|---|---|---|---|---|---|
| 10 | Adam optimizer | tanh | PI- | 0.000118 | 0.000125 |
| | | | PI+ | 0.000128 | 0.000144 |
| | Adam optimizer | ReLU | PI- | 0.000519 | 0.000560 |
| | | | PI+ | 0.000456 | 0.000496 |
| | CONOPT | tanh | PI- | 0.000020 | 0.000021 |
| | | | PI+ | 0.000018 | 0.000018 |
| 15 | Adam optimizer | tanh | PI- | 0.000054 | 0.000061 |
| | | | PI+ | 0.000063 | 0.000067 |
| 20 | Adam optimizer | tanh | PI- | 0.000047 | 0.000051 |
| | | | PI+ | 0.000048 | 0.000051 |

Table 7 summarizes the optimization performance for this case study, where different number of neurons, training algorithms and conditions, activation functions and optimization methods are used.

**Table 7.**

Optimization results for case study 3.

| *# of neurons in the hidden layer* | *Training algorithm* | *Objective function for training* | *Hidden layer activation function at training* | *Optimization method* | *# of feasible runs out of 5* | *Deviation ratio of $F_{RSD}$ from the global solution* | *Average SSE ($x10^{-4}$)* |
|---|---|---|---|---|---|---|---|
| 10 | Adam optimizer | (Eq.2) (PI-) | tanh | 3-pcs tanh | 4 | 0.08 | 5.3 |
| | | (Eq.3) (PI+) | | | 4 | 0.06 | 4.8 |
| | | (Eq.2) (PI-) | | 5-pcs tanh | 3 | 0.05 | 4.4 |
| | | (Eq.3) (PI+) | | | 4 | 0.04 | 4.0 |
| | | (Eq.2) (PI-) | ReLU | relumip | 5 | 0.00 | 9.3 |
| | | (Eq.3) (PI+) | | | 5 | 0.00 | 6.6 |
| | CONOPT | (Eq.2) (PI-) | tanh | 3-pcs tanh | 3 | 0.13 | 38.2 |
| | | (Eq.5) (PI+) | | | 5 | 0.07 | 34.7 |
| 15 | Adam optimizer | (Eq.2) (PI-) | tanh | 3-pcs tanh | 4 | 0.14 | 7.7 |
| | | (Eq.3) (PI+) | | | 4 | 0.08 | 3.2 |
| | | (Eq.2) (PI-) | | 5-pcs tanh | 4 | 0.07 | 3.9 |
| | | (Eq.3) (PI+) | | | 4 | 0.04 | 2.2 |
| 20 | Adam optimizer | (Eq.2) (PI-) | tanh | 3-pcs tanh | 4 | 0.16 | 7.9 |
| | | (Eq.3) (PI+) | | | 5 | 0.08 | 6.3 |
| | | (Eq.2) (PI-) | | 5-pcs tanh | 4 | 0.10 | 3.9 |
| | | (Eq.3) (PI+) | | | 5 | 0.04 | 4.3 |

### 3.3.1. Bi-objective physics-informed training and piecewise linear optimization

Built-in custom loss function in Keras library, which is used in this case study, computes predicted and true values (Chollet et al, 2015). Physics term is the square root of the differences of the flowrates of the inputs and outputs of the distillation column, which indicates the square root of mass balance:

$$p = \sqrt{\sum_{i=1}^{N} u_{1,i} - \sum_{j=1}^{6} \hat{y}_{i,j}} \tag{17}$$

where $u_{1,i}$ is the first ANN input of the i$^{th}$ data, representing the crude oil flowrate, and $\hat{y}_{i,j}$ is the j$^{th}$ column of i$^{th}$ data of the ANN outputs, representing the flowrate of the white products obtained from the distillation process.

The neural network models are trained in Tensorflow with a single hidden layer, using 10, 15 and 20 neurons. Adam optimizer is used as the training algorithm (Kingma & Ba, 2014). For 10, 15 and 20 neurons, hyperbolic tangent is used as activation function at the hidden layer. 5 runs have been performed for each case, and average MSE values for training and test are reported in Table 6.

For each number of neurons in the hidden layer, MSE values of training may slightly change in case of physics informed training. Weight and bias parameters obtained from training part are given as parameters to the optimization where the objective is to minimize the flowrate of residuum. This way, the optimal inputs that achieve this objective are computed. For the optimization, 3 and 5-segment piecewise linear form of hyperbolic tangent function is used for the approximation to MILP for the trained models with tanh as activation function.

In this case study, problem size is larger with higher number of neurons, which causes the optimization problem formulated with exact tanh quite challenging to solve (Schweidtmann & Mitsos, 2019). Corresponding MILP problems with piece-wise linear neural networks embedded are solved by the CPLEX solver in GAMS. The optimal input set values are then utilized in the mathematical unit model obtained from (Alhajri et al., 2008) to obtain the optimal outputs, and these outputs are expressed as first-principles model outputs ($y_{FP}$). At the final step, if the optimal set of the ANN surrogate model is feasible, the optimal outputs computed with ANN model optimization ($y^*$) and first-principles model results corresponding to the optimal input set obtained from piece-wise linear optimization, are compared by calculating the sum of squared errors (SSE), given in Eq. 18. This comparison criterion is important when the impact of suboptimal training is considered. The calculation scheme is given in Fig. 5.

$$SSE = \sum_{i=1}^{N}\sum_{j=1}^{6}\left(y_{FP,ij} - y_{i,j}^{*}\right)^{2} \tag{18}$$

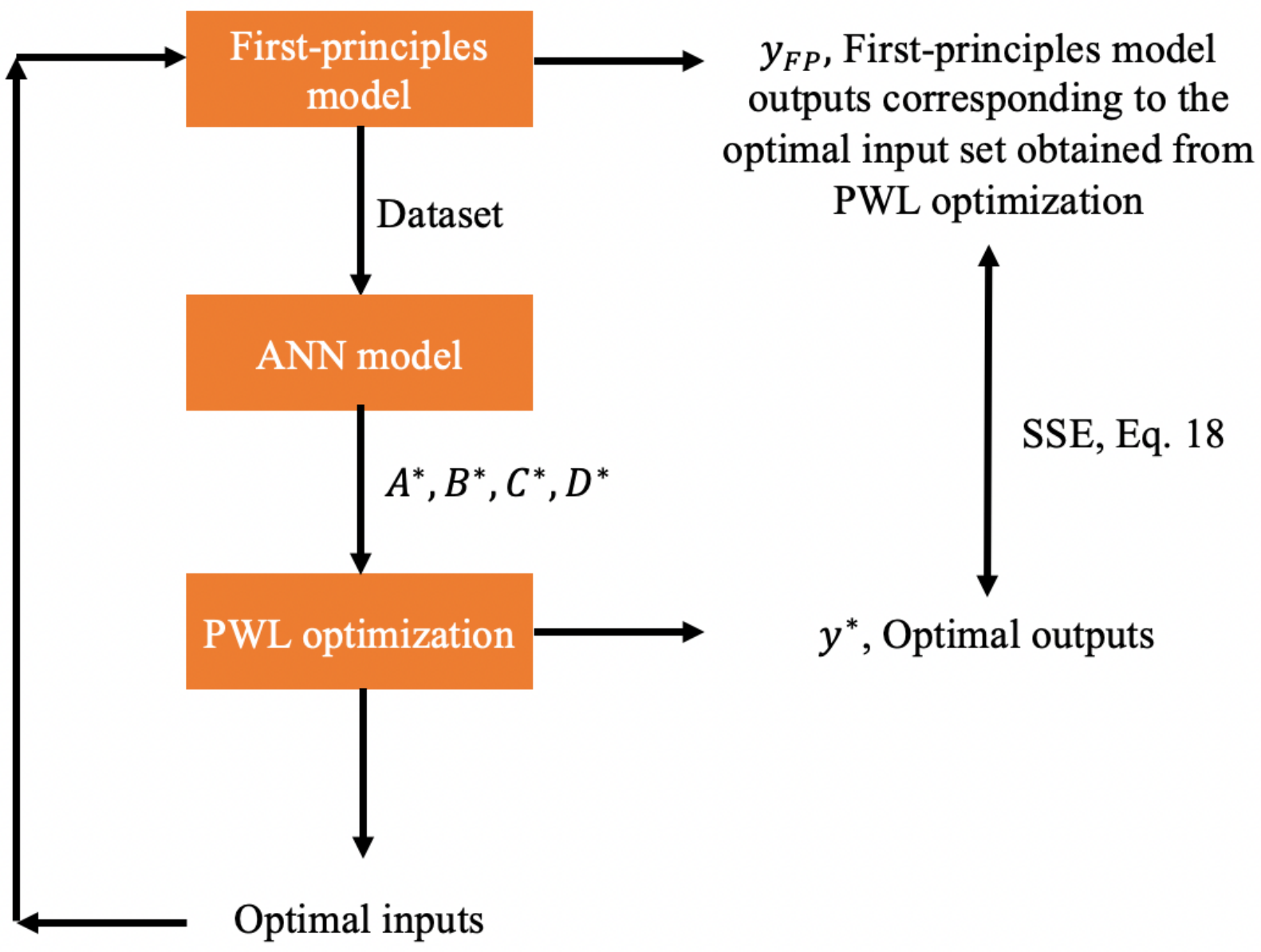

**Fig. 5.** Comparison scheme.

For different number of neurons in the hidden layer, average SSE value of feasible runs is taken for both physics informed and uninformed cases and are given in Table 7.

Regarding the models trained with tanh as activation function and Adam optimizer, there are three main outcomes:

(i) For the same number of neurons, physics-informed case yields closer results to the global solution. Although training and test MSE values are extremely closer, optimization results might change, showing that robust and consistent training may not even be enough for accurate surrogate optimization. It is observed that physics-informed training enables more reliable optimization.

(ii) ANN surrogate model optimization outputs are closer to the first-principles model outputs that correspond to the same optimal input set, when 5-piece form of tanh is used. This points out that instead of formulating the optimization problem with exact tanh to be solved as nonlinear programming (NLP) problem, which is quite challenging especially for higher number of neurons and bigger problem size, MILP formulations can be extended with larger number of breakpoints and segments, taking the advantage of the convexity.

(iii) When the optimal input set obtained from the surrogate optimization is utilized by the first-principles mathematical model, infeasibility may be observed, which is more probable when the neural network is trained with a physics-uninformed objective. This points out that, although the solutions obtained from the optimization are similar for PI- and PI+ cases, there is a possibility of obtaining physically unrealizable input set for PI- training and optimization. Another outcome is that physics-informed objective may overcome the impact of suboptimal training.

### 3.3.2. Multi-objective training and ReLU activation function

In this part, a benchmark between the standard usage of ReLU and newly proposed tanh approximations is provided. As mentioned earlier, using the ReLU activation function is the state-of-the art for such surrogate optimization studies. In this part, training is employed using 10 neurons in the hidden layer, where ReLU is used as the activation function and output layer utilizes the identity function. MSE values listed in Table 6 shows that training and test performances of physics informed case are better compared to the uninformed one when ReLU is used.

The reluMIP toolbox can be used for the surrogate optimization with neural network models embedded (Lueg et al., 2021). This toolbox uses the advantage of the exact piecewise nature of ReLU and formulates the optimization problem as an MILP. ReLU activation function is the mere choice though. By using Gurobi as MILP solver, optimal inputs that minimizes the crude oil flowrate are computed. Then, the optimal inputs are given to the first-principles model to obtain the outputs corresponding to the optimal input set. Average SSE values in Table 7 show that physics-informed training yields much less deviations and thus more accurate optimization results. Accordingly, it can be stated that physics-informed training is not only beneficial for tanh activation functions but also for ReLU, where publicly available open-source toolboxes can be used for optimization.

### 3.3.3. Constrained training

The ANN model is also trained in GAMS, where constrained training can be performed for PI+ case, addressed as physics-constrained training. For both PI- and PI+ cases, there are 10 neurons in the hidden layer where hyperbolic tangent is used as the activation function. CONOPT is used as local training solver, since traditional Tensorflow training algorithms are typically unable to solve constrained nonlinear programming problems.

PI- case is trained according to Eq. 2, where there is no physics term or constrain, whereas PI+ case is trained in two steps as follows:

a) First, there is no physics-term or constrain in the training, as in Eq. 5a, which corresponds to a warm-starting. This part is the same as physics-uninformed training.
b) In the second part, MSE and physics term constrains, $Z$ and $U$, are introduced as upper bounds, as in Eq. 5b. Weight and bias parameters are calculated again based on the MSE and physics constrains.

Table 6 shows that training and test instances of the physics-informed case perform better in terms of mean squared error compared to the physics-uninformed cases. This trend is different from other training and test performances of the neural networks trained with hyperbolic tangent at hidden layer, due to the advantages of warm-starting and constrained training. The warm-starting in the first step sustains better initialization for the physics-informed training part. In the second part, weight and bias terms are regulated based on the physics constrain, which increases the performance of training and test. When an upper bound is defined for the physics-term, the training and test performance may decrease, resulting in underfitting. Generally, this is a typical issue for multi-objective optimization. To avoid this issue, it is crucial to introduce constrains for both physics term and MSE in the second part. Deciding on the upper bounds is another issue for constrained training since the training problem with tight bounds can converge to an infeasible solution. On the other hand, when bounds are relaxed, constrained training may have negligible impact. The upper bounds have to be decided by trial-and-error, until satisfactory training performance is obtained.

The optimization process is the same as in Section 3.3.1 and 3-piece approximation to tanh is used. Physics-informed case by constrained training yields more reliable optimization with lower deviations from the global solution, as well as physics-informed bi-objective training.

As training algorithms, CONOPT and Adam optimizer uses generalized reduced gradient and extended version of stochastic gradient descent, respectively. Compared to Adam optimizer, higher training and test accuracy is obtained with CONOPT, but the optimal results deviate from the first-principles model solution when the models are trained with CONOPT, which may be due to the difference in the training algorithms and instances. Still, it is yielded that the crucial point is that physics-informed trained neural networks yield higher accuracy and precision for surrogate optimization regardless of the training algorithm and type of the physics-informed term implementation (multi-objective or constrained).

## 4. Conclusion

Physics informed training for neural networks is a promising way to overcome the complexity of physics-based models and to improve the black-box nature of data-driven models. For three different case studies consisting of a simple blending process, an industrial crude oil distillation column and a crude oil unit, the impact of physics-informed training compared to physics-uniformed one on process optimization is investigated. For the first two case studies, test performance increases since the mass balance eliminates the impact of the noise in the data. When the neural network is physics-informed trained, process optimization that is formulated by exact hyperbolic tangent and solved by nonconvex global solver converges to global optimums. MILP formulations with 5-piece approximation to hyperbolic tangent activation function yield significantly closer results to formulations with exact tanh, showing the generalization capacity of the proposed methods for process optimization. In fact, global solution is obtained almost every time with the physics informed piece-wise linear neural networks for surrogate optimization. Furthermore, MILP formulation requires much less CPU time due to convexity.

For the crude oil unit, which is the third case study, different training and optimization methods based on solvers, nonlinear activation functions and piece-wise linear approximations are tested. For a single hidden layer and 10, 15, 20 neurons, Adam optimizer is used for training with the nonlinear activation function being as tanh. It is demonstrated that physics-informed cases result in more reliable optimization results for each number of neurons. 5-piece approximation to tanh converges to closer results to global optimum compared to the 3-piece approximation. The outcome is the same when ReLU is selected as the activation function at the hidden layer.

Finally, the presence of physics-term as a constraint or a secondary objective in the training problem boosts the performance of the optimization, eliminating the challenges based on the quality and number of data, which is an issue in machine learning applications. Furthermore, in addition to using naturally piece-wise linear activation function such as ReLU, hyperbolic tangent can also be a nice candidate for formulating MILP optimization problems with physics-informed neural networks embedded. Therefore, number of pieces in the approximation to hyperbolic tangent can be increased to obtain closer results to the optimization with the exact form of the hyperbolic tangent. Finally, using physics-informed piece-wise linear networks for process optimization significantly reduces the computational times, which is potentially an advantage for deep-neural networks and real time optimization and control.

**Acknowledgements**

We gratefully acknowledge TÜPRAŞ refinery and TÜPRAŞ R&D department for sharing their actual process data.